\documentclass[a4paper]{article} \usepackage{theorem}
\usepackage{amsfonts} 
\usepackage{times}
\usepackage{mathptm}
\newcommand{\R}{\mathbb{R}} 
\title{Oriented Lagrangian Orthogonal Matroid Representations}
\date{January 2000 (Revised September 2000)}
\author{Richard F. Booth}
\newcommand{\qed}{\hfill $\diamond $ \\ \vskip 5pt}
\newcommand{\defterm}[1]{\emph{#1}} \theorembodyfont{\slshape}
\newtheorem{lem}{Lemma} \theorembodyfont{\slshape}
\newtheorem{thm}[lem]{Theorem} 
\theorembodyfont{\upshape} \newtheorem{defn}{Definition}
\theorembodyfont{\slshape} 
\theorembodyfont{\slshape} \newtheorem{ax}{Axiom}
 \newcommand{\B}{\ensuremath{{\cal B}}}
 \newcommand{\djunion}{\sqcup}
\newcommand{\sgn}{{\rm sign}\ } \newcommand{\Pf}{{\rm Pf}}
\newcommand{\pf}{{\rm p}} 
\newcommand{\mod}{{\rm mod}}
\newcommand{\af}[2]{\ensuremath{a_{f_{#1}f_{#2}}}}
\newcommand{\cf}[2]{\ensuremath{c_{f_{#1}f_{#2}}}}
\newenvironment{proof}[1][Proof]{\paragraph{#1}}{\qed}
\newcommand{\mtx}[1]{{\bf #1}}
\newcommand{\showlabel}[1]{\label{#1}}

\begin{document}
\maketitle

Several attempts have been made to extend the theory of matroids (here
referred to as \defterm{ordinary} or \defterm{classical} matroids) to
theories of more general objects, in particular the \defterm{Coxeter
matroids} of Borovik, Gelfand and White (\cite{book}, first introduced
as $WP$-matroids in \cite{GS}), and the \defterm{$\Delta $-matroids}
and (equivalent but for notation) \defterm{symmetric matroids} of
Bouchet (see, for example, \cite{bou}). The special cases of Coxeter
matroids for the Coxeter groups $BC_n $ and $D_n $ and a maximal
parabolic subgroup are called \defterm{symplectic} and
\defterm{orthogonal} matroids respectively, and may be viewed as
collections of $k$-element subsets of the $2n$-element set $\{1,
\ldots, n, 1^*, \ldots, n^* \}$ with maximality conditions, where $k$
is between 1 and $n$. In the case where $k = n $, these structures are
called \defterm{Lagrangian} matroids and are isomorphic in a natural
way to Bouchet's symmetric matroids \cite{bgw, orth}, with orthogonal
matroids giving \defterm{even} symmetric matroids. Classical matroids
now appear as a special case of even Lagrangian matroids. A concept of
representation of even $\Delta $- and symmetric matroids by
skew-symmetric $n \times n$ matrices was developed in \cite{bou3}. In
turn, symplectic and orthogonal matroids may be represented by
$k$-dimensional totally isotropic subspaces of $2n$-dimensional
symplectic and orthogonal vector spaces \cite{bgw, orth}; it is from
this that the names of these structures arise.

Attempts have also been made to extend the (classical) theory of
oriented matroids to this larger concept. A theory of orientation of
Lagrangian symplectic matroids was presented in \cite{bbgw}. However,
in the case when the matroid is even (as all orthogonal matroids are),
this theory is both uninteresting and trivial; in particular, it is
uninteresting for classical matroids. In \cite{wpp}, Wenzel presents
an orientation concept for even $\Delta $-matroids, and their
representations, which includes classical oriented matroids as a
special case. In this paper we extend this theory to Lagrangian
orthogonal matroids and their representations, and give a completely
natural transformation from a representation of a classical oriented
matroid to a representation of the same oriented matroid embedded as a
Lagrangian orthogonal matroid. We are interested in representations of
Lagrangian matroids as isotropic subspaces because such
representations arise in the study of maps on surfaces \cite{maps,
cohom}, and also because of their natural connections with Schubert
cells. Since classical represented matroids correspond to thin
Schubert cells in the Grassmannian \cite{bor-gel:schu}, oriented
matroids provide a stratification of the Grassmannian finer than thin
Schubert cells but coarser than their connected components. Similarly,
these other concepts of orientation provide stratifications of
Lagrangian varieties which split thin Schubert cells into unions of
connected components.

\section{Matroids and Representations}
\showlabel{sec:matandrep}

In this section we recall definitions of classical, symmetric and
$\Delta $-matroids. We then briefly discuss representations of these
objects and connections between them, and give an alternative
definition of Lagrangian orthogonal matroids.

\subsection{Matroids}

Let $I=\{1, \ldots, n \}$. Let $I_k = \{A \subseteq I \mid \#A = k\}$,
the collection of $k$-element subsets of $I$ (we use
\defterm{collection} for a set of sets to avoid confusion). Set also
$I^* = \{1^*, \ldots, n^*\} $, and $J = I \djunion I^* $. We define
the involution $*$ on $J$ by setting $(i^*)^* = i$ for $i^* \in I^* $
and extend it to sets in the obvious way. A set $A \subseteq J$ is
said to be \defterm{admissible} if $A \cap A^* = \emptyset $, and we
set $J_k $ to be the collection of admissible $k$-subsets of $J$. The
\defterm{symmetric difference} of two sets $A$ and $B$ is written and
defined by
\[ A \Delta B = (A \cup B) \setminus (A \cap B). \]

Then $M$ is a (classical) matroid if and only if it satisfies Axiom
\ref{ax:clmatex} below.
\begin{ax}[Classical Basis Exchange] \showlabel{ax:clmatex}
  For $A, B \in \B$ and $i \in A \setminus B $, there exists $j \in B
  \setminus A $ such that $(A \Delta \{i, j\}) \in \B$.
\end{ax}

Bouchet, in \cite{bou}, defines a \defterm{$\Delta $-matroid} as a
collection \B\  of 
subsets of $I$, not necessarily equicardinal, satisfying the
following:
\begin{ax}[Symmetric Exchange Axiom] \showlabel{ax:delexax}
For $A, B \in \B$ and $i \in A \Delta B $, there exists $j \in B
\Delta A $ such that $(A \Delta \{i, j\}) \in \B$.
\end{ax}
It is thus immediately apparent that a classical matroid is also a
$\Delta$-matroid. Bouchet  goes on to define a \defterm{symmetric matroid}
as essentially a $\Delta $-matroid with bases extended to $n$ elements
by adding to $B \in \B$ all starred elements which do not appear,
unstarred, in $B$. Thus a symmetric matroid is a set $\B \subseteq J_n
$ satisfying:
\begin{ax} \showlabel{ax:symexax}
For $A, B \in \B$ and $i \in A \Delta B $, there exists $j \in B
\Delta A $ such that $(A \Delta \{i, j, i^* , j^* \}) \in \B$.
\end{ax}
We shall refer to these two axioms interchangeably as `the
symmetric exchange axiom' depending on the structure to which we
refer.

We shall now define classical, symplectic and orthogonal matroids in
terms of maximality properties.  These definitions are drawn from
\cite{orth}; equivalences with other popular definitions may also
be found there.  Recall that, given a partial ordering $\prec$ on a
set $X$, the \defterm{Gale ordering} on the set of $k$-element subsets
$X_k$ of $X$ is defined as follows: for $A, B \in J_k$, write
\[ A=\{a_1,a_2,\ldots, a_k\}, B=\{b_1,b_2,\ldots,b_k\}, \]
with $a_i\prec a_{i+1}$ and $b_i\prec b_{i+1}$ for $1\le i < k$.  Then
we write $A\prec B$ if $a_i\prec b_i$ for $1\le i\le k$.

By a \defterm{$B_n$-admissible ordering}, we mean a total ordering on $J$
satisfying $i\prec j$ if and only if $j^*\prec i^*$; that is, an
ordering of the form
\[ a_1 \prec a_2 \prec \cdots \prec a_n \prec 
   a_n^* \prec a_{n-1}^* \prec \cdots \prec a_1^*   \]
where $\{a_1,\ldots,a_n\}\subset J$ is an admissible set.
By a \defterm{$D_n$-admissible ordering}, we mean a partial ordering
on $J$ of the form 
\[ a_1 \prec a_2 \prec \cdots \prec 
\begin{array}{c} a_n \\  a_n^* \end{array}
   \prec a_{n-1}^* \prec \cdots \prec a_1^*   \]
where $\{a_1,\ldots,a_n\}\subset J$ is an admissible set. Now we have
the following (standard) definitions:
\begin{enumerate}
\item A collection $\B\subseteq I_k$ is a (classical) matroid if and
only if for every linear ordering $\prec$ of $I$ there exists some
$B\in \B$ such that $A\prec B$ for every $A\in \B$. 
\item A collection $\B\subseteq J_k$ is a \defterm{symplectic matroid}
if and only if for every $B_n$-admissible ordering $\prec$ of $J$
there exists some $B\in \B$ such that $A\prec B$ for every $A\in \B$.
\item A collection $\B\subseteq J_k$ is an \defterm{orthogonal
matroid} if and only if for every $D_n$-admissible ordering
$\prec$ of $J$ there exists some $B\in \B$ such that $A\prec B$ for
every $A\in \B$. 
\end{enumerate}
Clearly, every orthogonal matroid is also a symplectic matroid.
A \defterm{Lagrangian} matroid is a symplectic matroid of maximal
rank (so that $k=n$).  Similarly, a Lagrangian orthogonal matroid is
an orthogonal matroid of maximal rank $n$, and Lagrangian orthogonal
matroids are Lagrangian matroids.

Finally, we observe that a Lagrangian (symplectic) matroid and a
symmetric matroid are the same objects.  This follows from the
characterisation of symmetric matroids in terms of a greedy algorithm
in \cite{bou}.  Furthermore, in \cite{orth}, it is shown that
Lagrangian matroids are orthogonal if and only if they are even; that
is, $B\cap I$ has the same parity for all bases $B$.  Thus, an
orthogonal Lagrangian matroid is exactly an even symmetric matroid.

\subsection{Representations}

\showlabel{sub:rep}

Concepts of representation of matroids have been introduced in two
separate, but closely related, ways. Bouchet introduces a concept of
representation by square matrices of `symmetric type' (\cite{bou3}),
whereas in \cite{bgw} representations are introduced in terms of
isotropic subspaces. In this paper we are concerned mainly with
representations over the real numbers.

Representable symplectic matroids arise naturally from symplectic and
orthogonal geometries, similarly to the way that classical matroids
arise from projective geometry. 

\paragraph{Classical representations}
We consider a $k$-dimensional subspace $U$ of a vector space $V$ with
basis $E=\{e_1, \ldots, e_n\}$. Choose a basis $u_1, \ldots, u_k$ for
$U$ and express it in terms of $E$ so that $u_i = \sum_{j=1}^n c_{i
j}e_j$. Thus, we have expressed this subspace as the row-space of a
$k\times n$ matrix $\mtx{C}$ of rank $k$ with columns indexed by $I$.
Let \B\ be the collection of sets of column indices corresponding to
non-zero $k\times k$ minors; then
\begin{thm} \showlabel{thm:repclsmat}
\B\  is the collection of bases of a (classical) matroid.
\end{thm}
Note that the matroid is independent of the choice of basis $u_1,
\ldots, u_k$. This theorem may be found in any book on matroid theory,
for example \cite{whitebook}. We now state the corresponding
result for symplectic and orthogonal matroids.

\paragraph{Symplectic and orthogonal representations}
Let $V$ be a vector space with basis
\[ E = \{ e_1, \ldots, e_n, e_{1^*}, , \ldots, e_{n^*} \}. \]
Let $\cdot$ be a bilinear form on $V$, with the symbol $\cdot$
often suppressed as usual, with
\begin{eqnarray*}
e_i e_{i^*} = 1 && \hbox{ for all } i \in I \\ e_i e_j = 0 && \hbox{
for all } i,j \in J \hbox{ with } i \ne j^*.
\end{eqnarray*}
\begin{defn}
The pair $(V, \cdot)$ is called a \defterm{symplectic space} for
$\cdot$ antisymmetric and an \defterm{orthogonal space} for $\cdot$
symmetric. If the vector space is of characteristic 2, it is
symplectic. A subspace $U$ of $V$ is called \defterm{totally isotropic}
if $\cdot$ restricted to $U$ is identically zero. A
\defterm{Lagrangian subspace} is a totally isotropic subspace of
maximal dimension (easily seen to be $n$).
\end{defn}

Choose a basis $u_1, \ldots, u_k$ of such a totally isotropic subspace
$U$ and represent this basis in terms of $E$, so that
\[ u_i = \sum_{j=1}^n \left( a_{i j} e_j + b_{i j} e_{j^*} \right) . \]
Now we have represented $U$ as the row space of a $k \times 2n$ matrix
$\mtx{C}= (\mtx{A}, \mtx{B})$ with columns indexed by $J$. Let \B\  be
the collection of sets of
column indices corresponding to non-zero $k\times k$ minors which are
admissible; then
\begin{thm} \showlabel{thm:repsymorthmat}
If $U$ is a totally isotropic subspace of a symplectic or orthogonal
space, \B\ is the collection of bases of a symplectic or orthogonal
matroid, respectively. Note that the matroid is independent of the
choice of basis $u_1, \ldots, u_k$ of $U$.
\end{thm}
This is Theorem 5 in \cite{orth}; the statement for symplectic
matroids only is Theorem 2 in \cite{bgw}.

$C$ is called a \defterm{(symplectic/orthogonal) representation} of
$M=(J, *, \B)$, and $M$ is said to be
\defterm{(symplecticly/orthogonally) representable}. Note that
orthogonal matroids may have symplectic representations. We also note
that, when considered in matrix form, the requirement that $U$ be
totally isotropic is equivalent to the requirement that $\mtx{AB^t}$
be symmetric in the symplectic case and skew-symmetric in the
orthogonal case.

In \cite{bou3}, Bouchet considers representations of $\Delta$-matroids
in terms of matrices of `symmetric type'.
\begin{defn} \showlabel{def:symtype}
A square matrix $\mtx{A}=(a_{i j})$ is said to be \defterm{quasi-symmetric}
if there exists a function $\epsilon : I \rightarrow \{-1, 1\}$ such
that $\epsilon(i) a_{i j} = \epsilon(j) a_{j i}$ for every $i,j \in
I$. Thus symmetric matrices are quasi-symmetric. $\mtx{A}$ is said to be of
\defterm{symmetric type} if it is anti-symmetric or quasi-symmetric.
\end{defn}
A \defterm{principal minor} of a square matrix is one consisting of
those rows and columns indexed by the same set $H\subseteq I$. Bouchet
proves
\begin{thm} \showlabel{thm:bourep}
Let the collection of subsets of $I$ corresponding to non-zero
principal minors of a matrix $\mtx{A}$ of symmetric type be $S$, and take any
$T\subseteq I$. Then the collection $\B = \{ A \Delta T \mid A\in S\}$
forms a $\Delta$-matroid.
\end{thm}
This is part of Theorem 4.1 in \cite{bou3}.

In fact, this result follows as a corollary of Theorem
\ref{thm:repsymorthmat}, and we can extend it a little in consequence.
Take a representation $\mtx{C}=(\mtx{A}, \mtx{B})$ of a Lagrangian
matroid $M$, choose a basis $F$ of it, and set $T = F \cap I$.
Exchange columns $j, j*$ for $j\in T$, and in the symplectic case
multiply one of each pair exchanged by $-1$. We have now moved those
columns corresponding to $F$ into the right-hand side while
maintaining (skew-) symmetry of $\mtx{AB^t}$. Now reduce, by row
operations, this non-singular right-hand-side to the identity matrix.
The resulting left-hand-side $\mtx{A}$ is clearly a symmetric matrix
in the symplectic case, and skew-symmetric in the orthogonal case.
This is now exactly the $\mtx{A}$ and $T$ of the above theorem. Other
sorts of quasi-symmetric matrices correspond to cases where the
right-hand-side has been reduced to a diagonal matrix with entries
plus or minus one, and indeed we may alter the definition of
`symmetric type' to read simply $\epsilon(i) a_{i j} = s\, \epsilon(j)
a_{j i}$, where $s=1$ or $s=-1$. We observe that any such
representation is equivalent to one which is strictly symmetric (for
$s=1$) or skew-symmetric (for $s=-1$) and that these produce
symplectic and orthogonal Lagrangian matroids respectively.

Note that we can `embed' a representation of a classical matroid as a
representation of the canonically associated Lagrangian orthogonal
matroid. (The classical matroid is a $\Delta$-matroid, which is a
symmetric matroid upon `completing' all sets in \B\  with the
appropriate starred elements. Since it is even, it is an orthogonal
Lagrangian matroid.) We simply make the top $k$ rows of $\mtx{A}$ (for
a matroid of rank $k$) the representation of the classical matroid,
and the remaining rows of $\mtx{A}$ zero; and the top $k$ rows of
$\mtx{B}$ zero, and the bottom $n-k$ rows an orthogonal complement of
maximal rank of $\mtx{A}$. This is clearly the required
representation, and is both a symplectic and an orthogonal
representation simultaneously.

In the case of a general, symplectically represented, symplectic
Lagrangian matroid, we assign orientations by considering essentially
signs of determinants of principal minors of the above symmetric matrices
\cite{bbgw}. Unfortunately, in skew-symmetric matrices that produces
uninteresting results, as we shall see; the correct concept is that of
the Pfaffian, which we shall define in the next section.

\section{Orientations}
\showlabel{sec:ori}

In this section we shall state a definition of classical oriented
matroids, give Wenzel's definition of (even) oriented
$\Delta$-matroids, and extend it in the obvious way to orthogonal
Lagrangian matroids. We remark parenthetically that symplectic
Lagrangian matroids (and so $\Delta$-matroids, even or otherwise) may
be oriented as described in \cite{bbgw}. We go on to discuss
representations of these objects, and prove that a representable
(classical) oriented matroid is representable as an oriented
orthogonal matroid.

\subsection{Orientation Axioms}
\showlabel{sub:oriax}

We begin by stating the Grassmann-Pl\"ucker relations.
\begin{thm} \showlabel{thm:grplrel}
  For all vectors
  $ x_1, \ldots, x_k, y_1, \ldots, y_k \in \R^k $
  we have that 
  \begin{eqnarray*} && \det (x_1, x_2, x_3, \ldots, x_k) \cdot \det (y_1, y_2, y_3,
    \ldots, y_k ) \\
    && \qquad =\sum_{i=1}^{k} \det (y_i, x_2, x_3, \ldots, x_k) \cdot \det
  (y_1, \ldots, y_{i-1}, x_1, y_{i+1}, \ldots, y_k ) \end{eqnarray*}
\end{thm}
The proof of this is simple: observe that the difference of the two
sides is an alternating multilinear form in the $k+1$ arguments $x_1,
y_1, y_2, \ldots, y_k$, vectors in a $k$-dimensional space. Hence this
form is zero.

These relations inspire the chirotope axioms of classical oriented
matroid theory:
\begin{defn} \showlabel{def:chirotope}
  A \defterm{chirotope} of rank $k$ on $I$ is a mapping $\chi: I^k
  \rightarrow  \{-1, 1, 0\}$ which satisfies:
  \begin{enumerate}
    \item $\chi$ is not identically zero.
    \item $\chi$ is alternating; that is
      \[ \chi (x_{\sigma(1)}, \ldots, x_{\sigma(k)}) = \sgn(\sigma) \chi (x_1, \ldots, x_k) \]
    for any $x_1, \ldots, x_k \in I,\ \sigma \in Sym(k)$.
    \item For all $x_1, \ldots, x_k, y_1, \ldots, y_k \in I$ such that 
      \[ \chi (y_i, x_2, x_3, \ldots, x_k) \cdot \chi (y_1, \ldots,
         y_{i-1}, x_1, y_{i+1}, \ldots, y_k ) \geq 0\]
      for $i=1, \ldots, k$ we have
      \[ \chi (x_1, x_2, x_3, \ldots, x_k) \cdot \chi (y_1, y_2, y_3,
    \ldots, y_k ) \geq 0. \]
  \end{enumerate}
\end{defn}
We then define an \defterm{oriented matroid} as an equivalence class
of chirotopes, where two chirotopes are said to be
\defterm{equivalent} if $\chi_1 = \pm \chi_2$. See \cite{oribook} for
a fuller description of this and other classical oriented matroid
definitions. We shall often speak of a chirotope as being an oriented
matroid, leaving the equivalence class implicitly understood.

We shall follow Wenzel in \cite{wpp} by making:
\begin{defn} \showlabel{def:twipf}
  A map $\pf: 2^I \rightarrow \R$ is called a \defterm{twisted Pfaffian map} if
  it satisfies the following:
  \begin{enumerate}
    \item $\pf$ is not identically zero.
    \item For all $A, B \subseteq I$ with $\pf(A)\ne 0$, $\pf(B) \ne
  0$, we have $\# A = \# B \ \mod \ 2$.
    \item If $A, B \subseteq I$ and $A\Delta B = \{i_1 < \ldots <
  i_l\}$ then we have
    \[ \sum_{j=1}^l (-1)^j \pf(A\Delta \{i_j\}) \cdot \pf(B\Delta
  \{i_j\}) = 0.\]
  \end{enumerate}
\end{defn}
We call two twisted Pfaffian maps \defterm{equivalent} if they differ
only by a non-zero constant scalar multiple. In fact, Wenzel makes the
definition for a `fuzzy ring' rather than for the real numbers, but we
are interested in this paper only in representations over the real
numbers. Pfaffian maps may be defined as twisted Pfaffian maps where
$\pf(\emptyset)=1$. 

\begin{defn}
Let
\[ S'_{2m} = \{ \sigma \in S_{2m} \mid \sigma(2k-1) = \min_{2k-1\le j
\le 2m} \sigma(j) \hbox{ for } 1\le k \le m\}, \]
and let $\mtx{A}$ be a skew-symmetric matrix.  Then the
\defterm{Pfaffian} of \mtx{A} is defined by 
\[
\Pf((a_{i j})_{1 \le i, j \le 2m}) = \sum_{\sigma \in S'_{2m}} \sgn
\sigma \prod_{k=1}^m a_{\sigma(2k-1)\,\sigma(2k)}.
\]
The Pfaffian of the empty set is 1, by definition.
\end{defn}
It can be shown that the square of the Pfaffian of a (skew-symmetric)
matrix is the determinant of that matrix.

\begin{thm} \showlabel{thm:wenzoriax} 
If $\mtx{A}$ is a skew-symmetric $n \times n$ matrix, $I_1, I_2 \subseteq I$
and $I_1 \Delta I_2 = \{ i_1,\ldots i_l \}$ with $i_j < i_{j+1}$ for $
1\le j \le l-1 $ then
$$
\sum_{j=1}^{l} (-1)^j \pf(I_1 \Delta \{i_j \}) \pf(I_2\Delta\{i_j
\})=0
$$
where $\pf(S) = \Pf((a_{i j})_{i, j \in S})$ for any $S \subseteq I$.
\end{thm}
This is Proposition 2.3 in \cite{w3}.

Thus a skew-symmetric matrix with real coefficients yields a Pfaffian
map, and in fact Pfaffian maps to a given ring (here, to the reals)
are in $1-1$ correspondence with skew-symmetric matrices over the same
ring (this is Theorem 2.2 in \cite{w3}). It is thus clear from Theorem
\ref{thm:bourep} that the subsets of $I$ corresponding to non-zero
values of the twisted Pfaffian map form a $\Delta$-matroid.

We now follow \cite[Definition 2.10]{wpp} in making
\begin{defn} \showlabel{def:oridelmat}
  An \defterm{oriented even $\Delta$-matroid} is an equivalence class of
  maps $\pf:2^I \rightarrow \{+1,-1,0\}$ satisfying
  \begin{enumerate}
    \item $\pf$ is not identically zero.
    \item For all $A, B \subseteq I$ with $\pf(A)\ne 0$, $\pf(B) \ne
  0$, we have $\# A = \# B\  \mod \ 2$.
    \item If $A, B \subseteq I$ and $A\Delta B = \{i_1 < \ldots <
  i_l\}$ and for some $w \in \{+1, -1\}$ we have
    \[ \kappa_j = w (-1)^j \pf(A\Delta \{i_j\}) \cdot \pf(B\Delta
  \{i_j\}) \geq 0\]
    for $1\leq j \leq l$, then $\kappa_j = 0$ for all $1\leq j\leq l$.
  \end{enumerate}
  We shall often speak of a map as an oriented even $\Delta$-matroid, with the
  equivalence class implicitly understood.
\end{defn}
The \defterm{bases} of the oriented even $\Delta$-matroid are those
subsets of $F\subseteq I$ for which $\pf(F) \neq 0$. We observe that
every Pfaffian map yields an oriented $\Delta$-matroid by simply
ignoring magnitudes.
\begin{lem}
The collection of bases of an oriented $\Delta$-matroid is a
$\Delta$-matroid.
\end{lem}
\begin{proof}
Recall that a collection of sets is a $\Delta$-matroid if and only if
it satisfies the symmetric exchange axiom, Axiom \ref{ax:delexax}:
\[ \hbox{for } E, F \in \B, e\in E\Delta F \hbox{ there exists } f\in
E\Delta F \hbox{ such that } E\Delta \{e,f\} \in \B.\] Set, without
loss of generality, $i_1=e$ in Condition 3 above (there is no loss of
generality since we are not concerned with signs or orderings). Set
$A=E\Delta \{e\}, B=F\Delta \{e\}$, and $w$ such that $\kappa_1$
is 1. Thus some other $\kappa_j$ must be $-1$; let $f=i_j$.  Now, from
the defining equation for $\kappa_j$, we have
\[ 0 \ne \pf(A\Delta \{f\})\pf(B\Delta \{f\}) =  \pf(E\Delta
\{e,f\})\pf(F\Delta \{e,f\}) \]
and so we obtain $E\Delta \{e,f\}, F\Delta \{e,f\} \in \B$,
which is more than we need.
\end{proof}
 We now make the obvious definition: Take a
Lagrangian orthogonal matroid \B, with an
equivalence class of signs assigned to its bases. Two sets of signs
are said to be equivalent when they are either identical on all bases
or opposite on all bases. We express this as an equivalence class of
maps
 $$\pf:J_n\rightarrow\{+,-,0\}$$
with $$\B = \{ A\in J_n \mid \pf (A) \neq 0\}$$ and equivalence given by $p\sim -p$.
Consider the
corresponding even $\Delta$-matroid and equivalence class of signs
$\pf'$ obtained by ignoring starred elements; that is, 
$\pf' (A) = \pf(B)$, where $B\in J_n$ is the unique element with
$B\cap I=A$. Now we say that $\pf$ is an \defterm{oriented orthogonal
matroid} exactly when $\pf'$ is an oriented even $\Delta$-matroid.

\subsection{Oriented Representations}
\showlabel{sub:orirep}

We first state two now-obvious theorems.
\begin{thm} \showlabel{thm:clsorirep}
  Given a $k \times n$ real matrix $\mtx{C}$, let 
\[ \chi (S \in I^k) = \sgn \det ( (c_{ij})_{j \in S} ). \]
  Then $\chi$ is an oriented matroid; further, the underlying
  (unoriented) matroid is the matroid represented by $\mtx{C}$. The oriented
  matroid represented is not altered when standard row operations are
  performed on $\mtx{C}$.
\end{thm}
\begin{thm} \showlabel{thm:oridelrep}
  Given an $n\times n$ square skew-symmetric real matrix $\mtx{A}$ and
  $T\subseteq I$, define $\pf : 2^I \rightarrow \{+1,-1,0\}$ by
  setting $\pf(B)$ to be the sign of the Pfaffian of the principal
  minor indexed by $B\Delta T$. Then $\pf$ is an oriented even
  $\Delta$-matroid, and the underlying $\Delta$-matroid is that
  represented by $\mtx{A}$ and $T$.
\end{thm}
The first theorem is classical, and the second from \cite{wpp}; both
should now be obvious from the definitions and earlier theorems.

An oriented classical matroid is described by a map
$$\chi:I^k\rightarrow \{+,-,0\}, \quad \chi \sim -\chi$$
 and an oriented even $\Delta$-matroid by a
map $$\pf:2^I\rightarrow \{+,-,0\}, \quad \pf \sim -\pf.$$ 
Given $\chi$, we widen the
domain by setting $\chi(A)=0$ whenever $\# A \neq k$, and obtain a map
which is a candidate to be an even $\Delta$-matroid. Given $\pf$
satisfying $\pf(A)=0$ whenever $\#A \neq k$, some fixed $k$, we can
restrict to a candidate to be an oriented matroid. It is natural to
ask when these candidates succeed.
\begin{thm} \showlabel{thm:clsoriisdel}
  Every oriented matroid is an oriented even $\Delta$-matroid, and
  every oriented even $\Delta$-matroid whose bases are all of rank $k$
  is an oriented matroid.

  Furthermore, a representation $\mtx{C}$ of an oriented matroid $M$
  yields a representation $\mtx{A}$ of it as an oriented even
  $\Delta$-matroid as follows. Choose a basis $T$ of $M$. Now set
  $a_{i j} = \det (T\Delta \{i,j\}) / \det(T)$, where by the
  determinant of a set we mean the determinant of the appropriate $k$
  columns of $\mtx{C}$, or $0$ if the set is not of cardinality $k$. Now
  $\mtx{A}$ is the required orientation.
\end{thm}
This follows from \cite[Theorem 4.1]{wpp}.

We now move on to define a representation of an oriented orthogonal
matroid.
\begin{defn} \showlabel{def:oriorthmat}
  Given $\mtx{C}$, an orthogonal representation of an orthogonal
  matroid $M$ over $\R$, we construct the oriented orthogonal matroid
  represented by $\mtx{C}$ as follows. Choose a basis $F$ of $M$, and
  swap columns $j$ and $j^*$ for $j\in T = F \cap I$ so that all
  columns of $F$ are in the right-hand $n$ places. Now perform row
  operations so that the right-hand $n$ columns become the identity
  matrix. Now the left-hand side, $\mtx{A'}$, is a skew-symmetric
  matrix (this is exactly the procedure discussed after Theorem
  \ref{thm:bourep}). Since we have $\mtx{A'}$ and $T$, we have a
  representation of an oriented even $\Delta$-matroid. Unfortunately,
  this oriented even $\Delta$-matroid is dependent on the initial
  choice of $F$, although the underlying non-oriented $\Delta$-matroid
  is not, so we modify $\mtx{A'}$ as follows.

  Set \[ \varepsilon_0 = 1 \hbox{ and } \varepsilon_i = \left\{
  \begin{array}{lr} \varepsilon_{i-1} & i \notin T \\
  -\varepsilon_{i-1} & i \in T \end{array} \right. \] for $i>0$. Then
  set $a_{ij} = \varepsilon_i \varepsilon_j a'_{ij}$.
  $\mtx{A}=(a_{ij})$ is again skew-symmetric, with rows and columns
  indexed by $I$, and we assign to the basis $B$ the sign of the
  Pfaffian of the principal minor of $\mtx{A}$ indexed by $(B\Delta
  F)\cap I$. If we consider instead that we have permuted column
  labels with columns, then the indices giving rise to this Pfaffian
  are those of the columns of $\mtx A$ labelled by elements of $B$.
  Note that this corresponds to the oriented even $\Delta$-matroid
  represented by $\mtx{A}, T$.
\end{defn}
\begin{thm} \showlabel{thm:oriorthmat}
  The above procedure obtains an oriented orthogonal matroid, which is
  independent of choice of $F$. 
\end{thm}
The fact that this is an oriented orthogonal matroid is obvious from
considering the oriented even $\Delta$-matroid represented by
$\mtx{A}, T$; we need only show independence of choice of $F$. It is
enough to show that a representation $(\mtx{A}, \mtx{I_n})$ yields the
same orientation using $\mtx{A}$ directly and going through the above
procedure with $\# T=2$. The symmetric exchange axioms of the first
section and the evenness tell us that any two bases are connected by a
path where adjacent bases differ in this way.

Suppose $a_{ij} \ne 0$, and set $T=\{i, j\}$ with $i<j$. Let $\mtx{B}$
be the skew-symmetric matrix obtained as follows. Take the compound
matrix $(\mtx{A}, \mtx{I_n})$, swap the $i$-th and $j$-th columns of
$\mtx{A}$ with those of $\mtx{I_n}$, and reduce using row operations
to the form $(\mtx{B},\mtx{I_n})$. It is helpful to know about the
form of $\mtx{B}$. When we write $\mtx{A}_S$, we mean the Pfaffian of
the minor of $\mtx{A}$ indexed by $S$. By $[a_{ij}a_{kl}]$, with $k\ne
l, \{i,j\}\cap\{k,l\}=\emptyset$, we mean $\pm \mtx{A}_{\{i,j,k,l\}}$
with the sign chosen such that the term $a_{ij}a_{kl}$ has positive
sign.
\begin{lem}\showlabel{lem:form}
The skew-symmetric matrix $B$ satisfies:
$$
  b_{kl}=\left\{ \begin{array}{lr}
  -1/a_{ij} & k=i, l=j, k\ne l \\
  a_{lj}/a_{ij} & k=i, l\ne j, k\ne l \\
  a_{il}/a_{ij} & k=j, l\ne i, k\ne l \\
  a_{jk}/a_{ij} & l=i, k\ne j, k\ne l \\
  a_{ki}/a_{ij}& l=j, k\ne i, k\ne l \\
  0 & k=l \\ 
  {[a_{ij}a_{kl}]}/a_{ij} &\quad |\{i,j,k,l\}|=4
  \end{array} \right.
$$
\end{lem}
\begin{proof}
 The first six statements are immediately clear from the construction
 of $\mtx{B}$. From consideration of determinants, which can be more
 readily seen, the final part is correct up to sign. But the term
 $a_{ij}a_{kl}$ appears in some sense `early' in the construction of
 $\mtx{B}$ from $\mtx{A}$ and cannot then change sign, so this is the
 correct sign also.
\end{proof}

Now, without loss of generality, $k<l$, since $b_{kl}=-b_{lk}$.
Define, for $|i,j,k,l|=4$, 
\[ \epsilon_{ijkl} = \left\{ \begin{array}{lr}
-1 & i<k<j<l \\ -1 & k<i<l<j \\ +1 & \hbox{ otherwise}
\end{array} \right.  \]
(we leave $\epsilon_{ijkl}$ undefined when its subscripts are not all
distinct). 
Clearly, from our formula for Pfaffians,
$[a_{ij}a_{kl}]=\epsilon_{ijkl}\mtx{A}_{\{i,j,k,l\}}$.

Let us define a matrix $\mtx{C}$ from $\mtx{B}$ by multiplying rows $k$ for $i\le
k < j$ and the corresponding columns by $-1$. Then we have
\begin{lem} \showlabel{lem:biglem}
 The Pfaffian minor $\mtx{C}_S$ satisfies
 $\mtx{C}_S=\mtx{A}_{S\Delta\{i,j\}}/a_{ij}$.
\end{lem}
\begin{proof}
 Throughout, $i<j$ and $k<l$. 
 Define $\rho_{ijk}=-1$ if
 $i\le k<j$ and $+1$ otherwise. Thus
 $\epsilon_{ijkl}=\rho_{ijk}\rho_{ijl}$ wherever $\epsilon_{ijkl}$ is
 defined, and $c_{kl}=\rho_{ijk}\rho_{ijl}b_{kl}$. Thus, the
 elements of the skew-symmetric matrix $\mtx{C}$ satisfy:
 \[ c_{ij}=1/a_{ij}\quad
   c_{ik}=\rho_{ijk}a_{jk}/a_{ij}\quad
   c_{jk}=\rho_{ijk}a_{ik}/a_{ij}\quad
   c_{kl} = \mtx{A}_{\{i,j,k,l\}}/a_{ij}. \]
   (for $|\{i,j,k,l\}|=4$).
 It is easy to see that the lemma holds for determinants rather than
 Pfaffians of minors, so $\mtx{C}_S=\pm \mtx{A}_S$.  Each term of
 $\mtx{C}_S$, rewritten in terms of the $a_{kl}$, corresponds to
 several terms of $\mtx{A}_S/a_{ij}$; thus we need check only that one
 of these has the same sign in $\mtx{C}_S$ as in $\mtx{A}_S$.

 Let
 $$f_1< \cdots < f_p < i < f_{p+1} < \cdots < f_q < j < f_{q+1} <
 \cdots < f_m$$  
 and write $F=\{f_1, \ldots, f_m\}$.
 We divide the proof into the four cases $S=F$, $S=F\cup\{i\}$,
 $S=F\cup\{j\}$, $S=F\cup\{i,j\}$.
 We shall divide these each into sub-cases depending on whether $p$
 and $q$ are odd or even.  

 First we take $S=F=\{f_1, \ldots, f_m\}$; we may assume $m$ is even (as
 otherwise $\mtx{A}_S=\mtx{C}_S=0$). Now, take
 $$\overline{c}=\cf{1}{2}\ldots\cf{m-1}{m},$$ which has positive sign
 in $\mtx{C}_S$; this contains the signed term
 $$\left(\prod_{t=1}^{m/2}
 \epsilon_{ijf_{2t-1}f_{2t}}\right)\overline{a}/a_{ij},$$ 
 where 
 \[\overline{a}=\af{1}{2}\ldots\af{m-1}{m}a_{ij}.\] 
 Now we consider our sub-cases. If both $p,q$ are even, then
 $\overline{a}$ has positive sign in $\mtx{A}_{S\Delta{i,j}}$, and in
 fact all the $\epsilon_{ijf_{2t-1}f_{2t}}$ are positive, so the term
 has positive sign in $\mtx{C}_S$ as well. If $p$ is odd but $q$ is
 even then $\overline{a}$ has negative sign in
 $\mtx{A}_{S\Delta{i,j}}$, and all the $\epsilon$ are positive except
 for $\epsilon_{ijf_{p}f_{p+1}}$, so again $\overline{c}$ has the
 correct sign. Similarly, if $q$ is odd but $p$ is even then
 $\overline{a}$ has negative sign, and all the $\epsilon$ are positive
 except for $\epsilon_{ijf_{q}f_{q+1}}$. Finally, if both $p,q$ are
 odd, then $\overline{a}$ has positive sign in
 $\mtx{A}_{S\Delta{i,j}}$, and all the $\epsilon$ are positive except
 for $\epsilon_{ijf_{p}f_{p+1}}$ and $\epsilon_{ijf_{q}f_{q+1}}$.
 This disposes of the first case.
 
 For the second case, take $S=F\cup\{i, j\}$.  Once again $m$ is even
 in the non-trivial case.
 Now $$\overline{a} = \af{1}{2}\ldots \af{m-1}{m}$$ has positive sign
 in $\mtx{A}_{S\setminus \{i,j\}}$, and
 $$\overline{c}=\cf{1}{2}\ldots\cf{m-1}{m}c_{ij}$$ yields the term
 $$\left(\prod_{t=1}^{m/2}
 \epsilon_{ijf_{2t-1}f_{2t}}\right)\overline{a}/a_{ij}.$$ Similarly to the
 first case, this
 is $\overline{a}/a_{ij}$ when $p,q$ are both even 
 or both odd, and $-\overline{a}/a_{ij}$ when exactly one of $p,q$ is even.
 However, $\overline{c}$ has positive sign in $\mtx{C}_S$ exactly when $p,q$ are
 both even or both odd. This disposes of the
 second case.
 
 Now take $S=F\cup\{i\}$.  Here the non-trivial case has $m$
 odd. Suppose first that $q$ is odd. Take
 $$\overline{a}=\af{1}{2}\ldots\af{q-2}{q-1}a_{f_q
 j}\af{q+1}{q+2}\ldots\af{m-1}{m},$$ 
 which has positive sign in $\mtx{A}_{F\cup \{j\}}.$
 Now take
 $$\overline{c}=\cf{1}{2}\ldots\cf{q-2}{q-1}c_{f_q
 i}\cf{q+1}{q+2}\ldots\cf{m-1}{m}.$$ This contains the term
 $$\left(\prod_{t=1}^{(q-1)/2} \epsilon_{ijf_{2t-1}f_{2t}}\right) \left(
 \prod_{t=(q+1)/2}^{(m-1)/2} \epsilon_{ijf_{2t}f_{2t+1}}\right)
 \rho_{ijf_q} \overline{a}/a_{ij}.$$ 
 Now, $\overline{c}$ has positive sign in
 $\mtx{C}_S$ exactly when $p$ is odd also. Since $i<f_q<j$,
 $\rho_{ijf_q}$ is negative, and all the $\epsilon$ are positive
 except for $\epsilon_{ijf_{p}f_{p+1}}$, which appears exactly when $p$ is
 odd.  This disposes of the sub-cases  where $q$ is odd.
 The remaining cases, for $q$ even and for $S=F\cup\{j\}$, are similar.
\end{proof}
Since $(\mtx{C}, \mtx{I_n})$ is the form that would be obtained by following
Definition \ref{def:oriorthmat}, we have proven Theorem
\ref{thm:oriorthmat}, as the signs differ only by a constant scalar multiple.
Finally, we state the following:
\begin{thm} \showlabel{thm:oriclstoorthrep}
  Let $\mtx{B}$ be a representation of the oriented matroid $M$, with
  columns indexed by $I$. Then \[ \left( \begin{array}{cc} \mtx{B} & 0
  \\ 0 & \mtx{D} \end{array} \right) \] is an orthogonal
  representation of the corresponding oriented orthogonal Lagrangian
  matroid, where $\mtx{D}$ is an orthogonal complement to $\mtx{B}$,
  with columns indexed by $I^*$.
\end{thm}
\begin{proof}
Let $M$ be of rank $k$, and suppose without loss of generality that
the leftmost $k$ columns of $\mtx{B}$ form a basis of $M$. Since performing
row operations on representations of classical oriented matroids does
not alter the oriented matroid represented, we may assume that these
$k$ columns form an identity matrix in the first $k$ rows, and that
the rightmost $n-k$ columns of the orthogonal complement form an
identity matrix in the last $n-k$ rows also. We swap these first $k$
columns into the right-hand-side, and make the appropriate
multiplications, obtaining a matrix $(\mtx A\ \mtx I)$, where
\[ \mtx A = \left( \begin{array}{cccccc}
0 &  \cdots & 0 & (-1)^k b_{1\ k+1}  & \cdots &  (-1)^k b_{1\ n} \\
\vdots  & \ddots & \vdots & \vdots \ddots & \vdots \\
0 & \cdots & 0 & (-1)^k b_{k\ k+1} & \cdots &  (-1)^k b_{k\ n} \\
(-1)^1 d_{1\ 1} & \cdots & (-1)^k d_{1\ k} & 0 & \cdots & 0 \\
\vdots & \ddots & \vdots & \vdots &  \ddots & \vdots \\
(-1)^1 d_{n-k\ 1} & \cdots & (-1)^k d_{n-k\ k} & 0 & \cdots & 0 \\
\end{array} \right). \]
Now we see that $a_{ij}=\det (\{1, \ldots, k\}
\Delta {i,j})$ where $\det$ is the determinant of the appropriate $k$
columns of $\mtx{B}$, and $0$ if its argument has more or less than $k$
elements. Now the result follows at once from Theorem \ref{thm:clsoriisdel}.
\end{proof}

\section*{Acknowledgement}

The author wishes to thank Neil White for his helpful advice, and for
proof-reading beyond the call of duty.

{\noindent
Richard F. Booth,
Department of Mathematics, UMIST, PO Box 88, Manchester M60 1QD,
United Kingdom; \tt richard.booth@umist.ac.uk}
\end{document}